\documentclass[12pt,a4paper]{amsart}
\usepackage[T1]{fontenc}
\usepackage[latin1]{inputenc}
\usepackage{geometry}
\makeatletter

\providecommand{\LyX}{L\kern-.1667em\lower.25em\hbox{Y}\kern-.125emX\@}

\theoremstyle{plain}    
\newtheorem{thm}{Theorem}[section]
\numberwithin{equation}{section} %% Comment out for sequentially-numbered
\numberwithin{figure}{section} %% Comment out for sequentially-numbered
\theoremstyle{plain}    
\newtheorem{cor}[thm]{Corollary} %%Delete [thm] to re-start numbering
\theoremstyle{plain}    
\newtheorem{lem}[thm]{Lemma} %%Delete [thm] to re-start numbering
\theoremstyle{plain}    
\newtheorem{prop}[thm]{Proposition} %%Delete [thm] to re-start numbering
\theoremstyle{remark}
\newtheorem{rem}[thm]{Remark}
\theoremstyle{remark}    
\newtheorem{notation}[thm]{Notation} 
\theoremstyle{remark}    
\newtheorem{acknowledgement}[thm]{Acknowledgement} 

\usepackage{amssymb}

\def\makebbb#1{
    \expandafter\gdef\csname#1\endcsname{
        \ensuremath{\Bbb{#1}}}
}
\makebbb{R}
\makebbb{N}
\makebbb{Z}
\makebbb{C}
\makebbb{G}
\makebbb{E}
\makebbb{H}
\makebbb{P}
\makeatother

\begin{document}

\title{Bergman Kernels and Local Holomorphic Morse Inequalities.}

\author{Robert Berman}

\email{robertb@math.chalmers.se}


\begin{abstract}
Let \( (X,\omega ) \) be a hermitian manifold and let \( L^{k} \) be a high
power of a hermitian holomorphic line bundle over \( X. \) Local versions of
Demailly's holomorphic Morse inequalities (that give bounds on the dimension
of the Dolbeault cohomology groups associated to \( L^{k}), \) are presented
- after integration they give the usual holomorphic Morse inequalities. The
local \emph{weak} inequalities hold on any hermitian manifold \( (X,\omega ) \),
regardless of compactness and completeness. The proofs, which are elementary,
are based on a new approach to pointwise Bergman kernel estimates, where the
kernels are estimated by a model kernel in \( \C ^{n}. \) 
\end{abstract}
\maketitle

\section{Introduction}

Let \( X \) be an \( n- \)dimensional (possibly non-compact) complex manifold
equipped with a hermitian metric two-form, denoted by \( \omega . \) Furthermore,
let \( L \) be a holomorphic line bundle over \( X. \) The hermitian fiber
metric on \( L \) will be denoted by \( \phi . \) In practice, \( \phi  \)
is considered as a collection of local functions. Namely, let \( s \) be a
local holomorphic trivializing section of \( L, \) then locally, \( \left| s(z)\right| _{\phi }^{2}=e^{-\phi (z)}, \)
and the canonical curvature two-form of \( L \) can be expressed as \( \partial \overline{\partial }\phi . \)\footnote{%
The normalized curvature two-form \( \frac{1}{\pi }\frac{i}{2}\partial \overline{\partial }\phi  \)
represents \( c^{1}(L), \) the first Chern class of \( L, \) in real cohomology.
} When \( X \) is compact Demailly's holomorphic Morse inequalities \cite{d1}
give asymptotic bounds on the dimension of the Dolbeault cohomology groups associated
to \( L^{k}: \) 
\begin{equation}
\label{demailly weak ineaq}
\dim _{\C }H_{\overline{\partial }}^{0,q}(X,L^{k}\otimes E)\leq k^{n}\frac{(-1)^{q}}{\pi ^{n}}\frac{1}{n!}\int _{X(q)}(\frac{i}{2}\partial \overline{\partial }\phi )^{n}+o(k^{n}),
\end{equation}
 where \( X(q) \) is the subset of \( X \) where the curvature-two form \( \partial \overline{\partial }\phi  \)
has exactly \( q \) negative eigenvalues and \( n-q \) positive ones. These
are the \emph{weak} holomorphic Morse inequalities - they also have \emph{strong}
counterparts involving alternating sums of the dimensions. Demailly's inspiration
came from Witten's analytical proof of the classical Morse inequalities for
the Betti numbers of a \emph{real} manifold \cite{wi}, where the role of the
fiber metric \( \phi  \) is played by a Morse function. Subsequently, proofs
based on asymptotic estimates of the heat kernel of the \( \overline{\partial }- \)
Laplacian, were given by Demailly, Bouche and Bismut (\cite{d2}, \cite{bo}
and \cite{bi}). All of these proofs use quite delicate analytical arguments
- heat kernel estimates and global estimates deduced from the Bochner-Kodaira-Nakano
identity for non-Kähler manifolds. In the present paper it is shown that Demaillys
inequalities, may be obtained from comparatively elementary considerations.
The starting point is the formula 
\begin{equation}
\label{dim B H as bergman kernel intro}
\dim _{\C }H_{\overline{\partial }}^{0,q}(X,L^{k})=\int _{X}\sum _{i}\left| \Psi _{i}(x)\right| ^{2},
\end{equation}
 where \( \left\{ \Psi _{i}\right\}  \) is any orthonormal base for the space
of \( \overline{\partial }- \)harmonic \( (0,q)- \)forms with values in \( L^{k}), \)
when \( X \) is compact (this is obvious for the dimension of the harmonic
space - and by the Hodge theorem the dimensions coincide). It is shown that
the integrand, called the \emph{Bergman kernel} \emph{function} \( B_{X}^{q,k}(x), \)
may be asymptotically estimated by a model kernel in \( \C ^{n}. \) Integration
then yields Demailly's weak inequalities and a similar argument gives the strong
inequalities. The main point of the proof is to first show the corresponding
localization property for the closely related \emph{extremal function \( S_{X}^{q,k}(x) \)}
defined as 
\[
\sup \frac{\left| \alpha (x)\right| ^{2}}{\left\Vert \alpha \right\Vert _{X}^{2}},\]
where the supremum is taken over all \( \overline{\partial }- \)harmonic \( (0,q)- \)forms
with values in \( L^{k}. \) Since the estimates are purely local, they hold
on any (possibly non-compact) complex manifold, and yield \emph{local} weak
holomorphic Morse inequalities for the corresponding \( L^{2}- \)objects. The
main inspiration for the present paper comes from Berntsson's recent article
\cite{be}. 

One final remark: it is fair to say that the formula \ref{dim B H as bergman kernel intro}
is the starting point for the previous writers approaches to Demailly's inequalities,
as well. The heat kernel approach is based on the observation that the term
corresponding to the zero eigen value in the heat kernel on the diagonal \( e^{q,k}(x,x;t) \)
is precisely the Bergman kernel function \( B^{q,k}_{X}(x) \) (if \( X \)
is compact). Moreover, when \( t \) tends to infinity the contribution of the
the other eigen values tends to zero. The main problem, then, is to obtain the
asymptotic expression for the heat kernel in \( k \) and \( t \) and investigate
the interchanging or the limits in \( k \) and \( t. \) The point of the present
paper is to work directly with the Bergman kernel.

\subsection{Statement of the main result\label{subsection statement of main result}}

Recall that the \( \overline{\partial }- \)Laplacian is defined by \( \Delta _{\overline{\partial }}:=\overline{\partial }\overline{\partial }^{*}+\overline{\partial }^{*}\overline{\partial } \)
(where \( \overline{\partial }^{*} \) denotes the formal adjoint of \( \overline{\partial }) \)
and its \( L^{2}- \)kernel is the space of harmonic \( (0,q) \)-forms with
values in \( L^{k}, \) denoted by \( \mathcal{H}^{0,q}(X,L^{k}). \) By the
well-known Hodge theorem this space is isomorphic to the Dolbeault cohomology
group \( H_{\overline{\partial }}^{0,q}(X,L^{k}) \), when \( X \) is compact.
Now define the \emph{Bergman kernel} \emph{function} \( B_{X}^{q,k}(x) \) and
the \emph{extremal function} \( S_{X}^{q,k}(x) \) of \( \mathcal{H}^{q}(X,L^{k}) \)
by
\[
B_{X}^{q,k}(x):=\sum _{i}\left| \Psi _{i}(x)\right| ^{2},\, \, \, \, \, S_{X}^{q,k}(x):=\sup \frac{\left| \alpha (x)\right| ^{2}}{\left\Vert \alpha \right\Vert _{X}^{2}},\]
 where \( \left\{ \Psi _{i}\right\}  \) is any orthonormal base in \( \mathcal{H}^{q}(X,L^{k}) \)
and the supremum is taken over all forms in \( \mathcal{H}^{q}(X,L^{k}). \)
To define the corresponding ``model'' functions at a point \( x \) in \( X, \)
\( B^{q}_{x,\C ^{n}} \) and \( S^{q}_{x,\C ^{n}}, \) proceed as before, replacing
the manifold \( X \) with \( \C ^{n}, \) the base metric in \( X \) with
the Euclidean metric in \( \C ^{n} \) and the fiber metric \( \phi  \) on
\( L \) with the fiber metric \( \phi _{0} \) on the trivial line bundle over
\( \C ^{n}, \) where
\[
\phi _{0}(z)=\sum ^{n}_{i=1}\lambda _{i,x}\left| z_{i}\right| ^{2},\]
and where the curvature two-form \( \partial \overline{\partial }\phi _{0} \)
in \( \C ^{n} \) is obtained by ``freezing'' the curvature two-form \( \partial \overline{\partial }\phi  \)
on \( X, \) at the point \( x \) (with respect to on orthonormal frame at
\( x). \)\footnote{%
Equivalently: the \( \lambda _{i,x} \) are the eigenvalues of the curvature
two-form of \( L \) with respect to on orthonormal frame at \( x. \) 
} Finally, denote by \( X(q) \) the subset of \( X \) where the curvature-two
form \( \partial \overline{\partial }\phi  \) has exactly \( q \) negative
eigenvalues and \( n-q \) positive ones. Its characteristic function is denoted
by \( 1_{X(q)}. \) Moreover, let \( X(\leq q):=\bigcup _{0\leq i\leq q}X(i). \)

The first theorem we shall prove is a local version of Demailly's weak holomorphic
Morse inequalities. Note that the manifold \( X \) is not assumed to be compact,
neither is there any assumption, e.g. completeness, on the hermitian metric
\( \omega  \). 

\begin{thm}
\label{thm:localweak}Let \( (X,\omega ) \) be a hermitian manifold. Then 

\[
\limsup _{k}\begin{array}{lr}
\frac{1}{k^{n}}B_{X}^{q,k}(x)\leq B^{q}_{x,\C ^{n}}(0), & \limsup _{k}\frac{1}{k^{n}}S_{X}^{q,k}(x)\leq S^{q}_{x,\C ^{n}}(0)
\end{array},\]
and 
\[
B^{q}_{x,\C ^{n}}(0)=S^{q}_{x,\C ^{n}}(0)=\frac{1}{\pi ^{n}}1_{X(q)}\left| det_{\omega }(\frac{i}{2}\partial \overline{\partial }\phi )_{x}\right| \]
 Moreover, 
\[
\lim _{k}\frac{1}{k^{n}}B_{X}^{q,k}(x)=\lim _{k}\frac{1}{k^{n}}S_{X}^{q,k}(x)\]
 if one of the limits exists.
\end{thm}
This seems to be a new result. The inequality for \( S_{X}^{q,k}(x) \) generalizes
results of Bouche \cite{bo} and Gillet, Soulé, Gromov \cite{gi}, that are
non-local and concern \emph{compact} manifolds \( X. \) Integration of the
inequality for the Bergman kernel gives Demailly's weak inequalities \ref{demailly weak ineaq}. 

When \( X \) is compact the local Morse inequalities can be extended to an
asymptotic \emph{equality} as follows. Let \( B_{\leq \nu _{k}}^{q,k} \) be
the Bergman kernel function of the space spanned by all the eigenforms of the
\( \overline{\partial }- \)Laplacian, whose eigenvalues are bounded by \( \nu _{k}. \)

\begin{thm}
Let \( (X,\omega ) \) be a compact hermitian manifold. Then 
\[
\lim _{k}k^{-n}B_{\leq \mu _{k}k}^{q,k}(x)=\frac{1}{\pi ^{n}}1_{X(q)}\left| det_{\omega }(\frac{i}{2}\partial \overline{\partial }\phi )_{x}\right| ,\]
for some sequence \( \mu _{k} \) tending to zero. 
\end{thm}
Again, integrating this yields, with the help of a well-known homological algebra
argument Demailly's \emph{strong} holomorphic Morse inequalities\footnote{%
Just as in the papers of Demailly, Bouche and Bismut (\cite{d2},\cite{bo}
and \cite{bi}). In fact, the idea of using the ``low-energy spectrum'' was
introduced in Witten's seminal paper \cite{wi}. A similar technique is used
in the heat equation proof of Atiyah-Singer's index theorem.
}. When the curvature of the line bundle \( L \) is strictly positive the asymptotic
equality holds for the usual Bergman kernel for the space of holomorphic sections
of \( L^{k}. \) This was first proved by Tian \cite{ti}) with a certain control
on the lower order terms in \( k. \) A complete asymptotic expansion was given
in \cite{z} using microlocal analysis. See also \cite{lin} where the manifold
is \( \C ^{n} \) and \cite{bor} where the complex structure is non-integrable.

\begin{rem}
The two theorems have straight forward generalizations to the case where the
forms take values in \( L^{k}\otimes E, \) for a given rank \( r \) hermitian
holomorphic vector bundle \( E \) over \( X. \) The estimate for the extremal
function \( S_{X}^{q,k} \) is unaltered, while the results for the Bergman
kernel function \( B_{X}^{q,k} \) are modified by a factor \( r \) in the
right hand side.
\end{rem}
\begin{notation}
The notation \( a\sim (\lesssim )b \) will stand for \( a=(\leq )C_{k}b, \)
where \( C_{k} \) tends to one with \( k. \)
\end{notation}

\subsection{A sketch of the proof of the local weak holomorphic Morse inequalities}

First we we will show how to obtain the estimate 
\begin{equation}
\label{weak inequality scetch}
\limsup _{k}k^{-n}S^{q,k}_{X}(x)\leq S^{q}_{x,\C ^{n}}(0)
\end{equation}
 By definition, there is a unit norm sequence \( \alpha _{k} \) of harmonic
forms with values in \( L^{k} \) such that 
\[
\limsup _{k}k^{-n}S_{X}^{q,k}(x)=\limsup _{k}k^{-n}\left| \alpha _{k}(x)\right| ^{2}\]
Now consider the restriction of the the form \( \alpha _{k}(x) \) to a ball
\( B_{R_{k}} \) with center in the point \( x \) and with radius \( R_{k} \)
decreasing to zero with \( k. \) The main point of the proof is that the form
\( \alpha _{k} \) is asymptotically harmonic, with respect to a \emph{model}
fiber metric, on a ball of radius slightly larger than \( \frac{1}{\sqrt{k}} \)
and is then a candidate for the corresponding model extremal function. Indeed,
we can arrange that the fiber metric \( k\phi  \) on the line bundle \( L^{k} \)
can be written as 
\[
\begin{array}{ccc}
k\phi (z)=k(\sum ^{n}_{i=1}\lambda _{i}\left| z_{i}\right| ^{2}+kO(\left| z\right| ^{3})) & 
\end{array}\]
in local coordinates around \( x. \) The form \( \beta ^{(k)}:=k^{-\frac{1}{2}n}\alpha ^{(k)}, \)
defined on the scaled ball \( B_{\sqrt{k}R_{k}}, \) where \( \alpha ^{(k)}(z) \)
denotes the (component wise) scaled form \( \alpha (\frac{z}{\sqrt{k}}), \)
satisfies
\[
\limsup _{k}k^{-n}S^{q,k}(x)=\limsup _{k}\left| \beta ^{(k)}(0)\right| ^{2}.\]
Moreover, \( \beta ^{(k)} \) is harmonic with respect to the scaled Laplacian
\( \Delta _{\overline{\partial }}^{(k)} \), taken with respect to the scaled
fiber metric \( (k\phi )^{(k)} \)on \( L^{k}: \) 
\[
(k\phi )^{(k)}(z)=\sum ^{n}_{i=1}\lambda _{i}\left| z_{i}\right| ^{2}+\frac{1}{\sqrt{k}}O(\left| z\right| ^{3})\]
 The point is that the scaled fiber metric converges to the quadratic model
fiber metric \( \phi _{0} \) in \( \C ^{n}, \) whith an appropriate choice
of the radii \( R_{k}. \) As a consequence the operator \( \Delta _{\overline{\partial }}^{(k)} \)
converges to the model Laplacian \( \Delta _{\overline{\partial },\phi _{0}}. \)
Standard techniques for elliptic operators then yields a subsequence of forms
\( \beta ^{(k_{j})} \) converging to a form \( \beta  \) defined in all of
\( \C ^{n}, \) which is harmonic with respect to the model Laplacian. This
means that

\[
\limsup _{k}k^{-n}S^{q,k}(x)=\left| \beta (0)\right| ^{2}\leq S^{q}_{x,\C ^{n}}(0),\]
proving \ref{weak inequality scetch}. Finally, lemma \ref{lemma: comparing B and S in general},
relating \( S^{q,k}(x) \) and the Bergman kernel function \( B^{q,k}(x), \)
is used to deduce the corresponding estimate for the Bergman kernel function.
All that remains is to compute the Bergman kernel and the extremal function
in the model case (section \ref{section:model}).

\section{\label{section: dimension, ...}The Bergman kernel function \protect\( B(x)\protect \)
and the extremal function \protect\( S(x)\protect \).}

In section \ref{subsection statement of main result} the Bergman kernel function
\( B_{X}^{q,k} \) and the extremal function \( S_{X}^{q,k} \) were defined.
We will also have use for component versions of \( S_{X}^{q,k}. \) For a given
orthonormal frame \( e_{x}^{I} \) in \( \bigwedge ^{0,q}_{x}(X,L^{k}) \) let
\[
S_{X,I}^{q,k}(x):=\sup \frac{\left| \alpha _{I}(x)\right| ^{2}}{\left\Vert \alpha \right\Vert _{X}^{2}},\]
where \( \alpha _{I}(x) \) denotes the component of \( \alpha  \) along \( e_{x}^{I}. \)
It will be clear from the context what frame is being used. These functions
are closely related according to the following, simple, yet very useful lemma.
Its statement generalizes a lemma used in \cite{be} (see also \cite{li}). 

\begin{lem}
\label{lemma: comparing B and S in general}Let \( L \) be a hermitian holomorphic
line bundle over \( X. \) With notation as above
\[
S_{X}^{q,k}(x)\leq B_{X}^{q,k}(x)\leq \sum _{I}S^{q,k}_{X,I}(x)\]

\end{lem}
\begin{proof}
To prove the first inequality in the statement, take any \( \alpha  \) in \( \mathcal{H}^{0,q}(X,L^{k}) \)
of unit norm. Since \( \alpha  \) is contained in an orthonormal base, obviously
\( \left| \alpha (x)\right| ^{2}\leq B(x), \) which proves the first inequality.
For the second inequality, let \( \left\{ \Psi ^{i}\right\}  \) be an orthonormal
base for \( \mathcal{H}^{0,q}(X,L^{k}) \), so that \( B_{X}^{q,k}(x):=\sum _{i}\left| \Psi ^{i}(x)\right| ^{2}=\sum _{i,J}\left| \Psi _{J}^{i}(x)\right| ^{2}. \)
Fix a \( J \) and let \( c_{i}:=\overline{\Psi _{J}^{i}(x).} \) Then, summing
only over \( i, \) gives 
\[
\sum _{i}\left| \Psi _{J}^{i}(x)\right| ^{2}=\sum _{i}c_{i}\Psi _{J}^{i}(x)=\alpha _{J}(x)\left( \sum _{i}\left| c_{i}\right| ^{2},\right) ^{1/2}\]
where \( \alpha =\sum _{i}\frac{c_{i}}{\left( \sum _{i}\left| c_{i}\right| ^{2},\right) ^{1/2}_{i}}\Psi ^{i} \)
lies in \( \mathcal{H}^{q}(X,L^{k}) \) and is of unit norm. Thus, 
\[
\left( \sum _{i}\left| \Psi _{J}^{i}(x)\right| ^{2}\right) ^{1/2}=\alpha _{J}(x)=\frac{\left| \alpha _{J}(x)\right| }{\left\Vert \alpha \right\Vert }\leq S_{J}(x)^{1/2}.\]
Finally, squaring the last relation and summing over \( J \) proves the second
inequality. 
\end{proof}

\section{\label{section:weakMorse}The Weak Holomorphic Morse Inequalities}

In this section we prove the local version of Demailly's holomorphic Morse inequalities
over \emph{any} complex manifold \( X. \) The usual version is obtained as
a corollary. Around each point \( x \) in \( X, \) fix local complex coordinates
\( \left\{ z_{i}\right\}  \) and a holomorphic trivializing section \( s \)
of \( L \) such that\footnote{%
this is always possible; see \cite{we}.
}
\[
\begin{array}{ccc}
\omega (z)=\frac{i}{2}\sum _{i,j}h_{ij}(z)dz_{i}\wedge \overline{dz_{j},}\textrm{ }\, h_{ij}(0)=\delta _{ij} & 
\end{array}\]
 and
\[
\begin{array}{ccc}
\left| s(z)\right| ^{2}=e^{-\phi (z)},\, \, \, \phi (z)=\sum ^{n}_{i=1}\lambda _{i,x}\left| z_{i}\right| ^{2}+O(\left| z\right| ^{3}), & 
\end{array}\]
where the quadratic part of \( \phi  \) is denoted by \( \phi _{0}. \) Note
that in the local coordinates the base metric \( \omega  \) coincides with
the Euclidean metric at the origin. Moreover, we have the identities,
\[
\lambda _{1,x}\lambda _{2,x}\cdot \cdot \cdot \lambda _{n,x}Vol_{(X,\omega )}=det_{\omega }(\frac{i}{2}\partial \overline{\partial }\phi )_{x}Vol_{(X,\omega )}=\frac{1}{n!}(\frac{i}{2}\partial \overline{\partial }\phi )^{n}.\]
 The following notation will be useful. Let \( B_{R}:=\{z:\, \left| z\right| <R\} \)
in \( \C ^{n} \) and let \( R_{k}:=\frac{\textrm{ln}k}{\sqrt{k}}. \) Using
the local chart around \( x, \) a (small) ball \( B_{R} \) is identified with
a subset of \( X \). Given a function \( f \) on the ball \( B_{R_{k}}, \)
we define a \emph{scaled} function on \( B_{\sqrt{k}R_{k}} \) by
\[
f^{(k)}(z):=f(\frac{z}{\sqrt{k}})\]
 Forms are scaled by scaling the components. Observe that scaling the fiber
metric on \( L^{k} \) gives 
\begin{equation}
\label{fiber metric expansion}
(k\phi )^{(k)}(z)=\sum ^{n}_{i=1}\lambda _{i}\left| z_{i}\right| ^{2}+\frac{1}{\sqrt{k}}O(\left| z\right| ^{3})
\end{equation}
 which motivates the choice of scaling.  The radius \( R_{k}:=\frac{\textrm{ln}k}{\sqrt{k}} \)
has been chosen to make sure that the fiber metric on \( L^{k} \) tends to
the quadratic model fiber metric \( \phi _{0} \) with all derivatives on scaled
balls, i.e.

\begin{equation}
\label{convergence of fiber m (new)}
\sup _{\left| z\right| \leq \sqrt{k}R_{k}}\left| \partial ^{\alpha }((k\phi )^{(k)}-\phi _{0})(z)\right| \rightarrow 0,
\end{equation}
where the convergence is of the order \( \frac{1}{\sqrt{k}} \) to some power,
which follows immediately from the expansion \ref{fiber metric expansion}.
Moreover, \( \sqrt{k}R_{k} \) tends to infinity, so that the sequence of scaled
balls \( B_{\sqrt{k}R_{k}} \) exhausts \( \C ^{n}. \) Let us denote by \( \Delta _{\overline{\partial }}^{(k)} \)
the Laplacian, taken with respect to the scaled fiber metric \( (k\phi )^{(k)} \)
and the scaled base metric \( \omega ^{(k)}. \) One can check that 
\begin{equation}
\label{scaling transf for laplace}
\Delta _{\overline{\partial }}^{(k)}\alpha =\frac{1}{k}(\Delta _{\overline{\partial }}\alpha )^{(k)}.
\end{equation}
Hence, if \( \alpha ^{k} \) is a form with values in \( L^{k}, \) which is
harmonic with respect to the global Laplacian, i.e \( \Delta _{\overline{\partial }}\alpha ^{k}=0, \)
then the scaled form \( \alpha ^{(k)} \) satifies 
\[
\Delta _{\overline{\partial }}^{(k)}\alpha ^{(k)}=0,\]
 on the scaled ball \( B_{\sqrt{k}R_{k}}. \) Moreover, because of the convergence
property \ref{convergence of fiber m (new)} it is not hard to check that 
\begin{equation}
\label{operator localization}
\Delta ^{(k)}_{\overline{\partial }}=\Delta _{\overline{\partial }}+\epsilon _{k}\mathcal{D}_{k},
\end{equation}
 where \( \mathcal{D}_{k} \) is a second order partial differential operator
with bounded variable coeffiecients on the scaled ball \( B_{\sqrt{k}R_{k}} \)
and \( \epsilon _{k} \) is a sequence tending to zero with \( k \) (in fact,
all the derivatives of the coefficients of \( \mathcal{D}_{k} \) are uniformly
bounded). It also follows from \ref{convergence of fiber m (new)} that for
any form \( \alpha ^{k} \) with values in \( L^{k}, \) 
\begin{equation}
\label{norm localization}
\left\Vert \alpha _{k}\right\Vert _{B_{R_{k}}}\sim k^{-n}\left\Vert \alpha ^{(k)}\right\Vert _{\phi _{0},\sqrt{k}R_{k}},
\end{equation}
 where the factor \( k^{-n} \) comes from the change of variables \( w=\frac{z}{\sqrt{k}} \)
in the integral. 

The proof of the following lemma is based on standard techniques for elliptic
operators.

\begin{lem}
\label{lemma:bootstrap}For each \( k, \) suppose that \( \beta ^{(k)} \)
is a smooth form on the ball \( B_{\sqrt{k}R_{k}} \) such that \( \Delta _{\overline{\partial }}^{(k)}\beta ^{(k)}=0. \)
Identify \( \beta ^{(k)} \) with a form in \( L_{\phi _{0}}^{2}(\C ^{n}) \)
by extending with zero. Then there is constant \( C \) independent of \( k \)
such that 

\[
\sup _{z\in B_{1}}\left| \beta ^{(k)}(z)\right| _{\phi _{0}}^{2}\leq C\left\Vert \beta ^{(k)}\right\Vert _{\phi _{0,}B_{2}}^{2}\]
Moreover, if the sequence of norms \( \left\Vert \alpha ^{(k)}\right\Vert _{\phi _{0},\C ^{n}}^{2} \)
is bounded, then there is a subsequence of \( \left\{ \beta ^{(k)}\right\}  \)
which converges uniformly with all derivatives on any ball in \( \C ^{n} \)
to a smooth form \( \beta , \) where \( \beta  \) is in \( L_{\phi _{0}}^{2}(\C ^{n}). \)
\end{lem}
\begin{proof}
Fix a ball \( B_{R}, \) of radius \( R \) in \( \C ^{n}. \) By Gårding's
inequality for the elliptic operator \( (\Delta _{\overline{\partial }}^{(k)})^{m}, \)
we have the following estimates for the Sobolev norm of \( \beta ^{(k)} \)on
the ball \( B_{R} \) with \( 2m \) derivatives in \( L^{2}: \) 
\begin{equation}
\label{proofLemmaGardingA}
\left\Vert \beta ^{(k)}\right\Vert _{B_{R},2m}^{2}\leq C_{R,k}\left( \left\Vert \beta ^{(k)}\right\Vert _{B_{2R}}^{2}+\left\Vert (\Delta _{\overline{\partial }}^{(k)})^{m}\beta ^{(k)}\right\Vert _{B_{2R}}^{2}\right) ,
\end{equation}
 for all positive integers \( m. \) Since \( \Delta _{\overline{\partial }}^{(k)} \)
converges to \( \Delta _{\overline{\partial },\phi _{0}} \) on the ball \( B_{2} \)
it is straight forward to see that \( C_{R,k} \) may be taken to be independent
of \( k. \) Hence, for all positive integers \( m, \) 
\[
\left\Vert \beta ^{(k)}\right\Vert _{B_{1},2m}^{2}\leq C\left\Vert \beta ^{(k)}\right\Vert _{B_{2}}^{2}\]
and the continuous injection \( L^{2,k}\hookrightarrow C^{0},\, k>n, \) provided
by the Sobolev embedding theorem, proves the first statement in the lemma. To
prove the second statement assume that \( \left\Vert \alpha ^{(k)}\right\Vert _{\phi _{0},\C ^{n}}^{2} \)
is uniformly bounded in \( k. \) Then \ref{proofLemmaGardingA} shows that
\[
\left\Vert \beta ^{(k)}\right\Vert _{B_{R},2,2m}^{2}\leq D_{R}\]
Since this holds for any \( m\geq 1, \) Rellich's compactness theorem yields,
for each \( R, \) a subsequence of \( \left\{ \beta ^{(k)}\right\}  \), which
converges in all Sobolev spaces \( L^{2,k}(B_{R}) \) for \( k\geq 0. \) The
compact embedding \( L^{2,k}\hookrightarrow C^{l},\, k>n+\frac{1}{2}l, \) shows
that the sequence converges in all \( C^{l}(B_{R}). \) Choosing a diagonal
sequence, with respect to a sequence of balls exhausting \( \C ^{n}, \) finishes
the proof of the lemma.
\end{proof}
Before turning to the proof of the local weak holomorphic Morse inequalities,
theorem, we state the following facts about the model case, that will be proved
in the following section:
\[
B^{q}_{x,\C ^{n}}(0)=S^{q}_{x,\C ^{n}}(0)=\frac{1}{\pi ^{n}}1_{X(q)}(x)\left| det_{\omega }(\frac{i}{2}\partial \overline{\partial }\phi )_{x}\right| .\]
 Moreover, suppose that the first \( q \) eigenvalues of the quadratic form
\( \phi _{0} \) are negative and the rest are positive (which corresponds to
the case when \( x \) is in \( X(q)). \) Then 
\begin{equation}
\label{S component in model}
S^{q}_{I,x,\C ^{n}}(0)=0,
\end{equation}
 unless \( I=(1,2,...,q). \)

\begin{thm}
Let \( (X,\omega ) \) be a hermitian manifold. Then the Bergman kernel function
\( B_{X}^{q,k} \) and the extremal function \( S_{X}^{q,k} \)of the space
of global \( \overline{\partial }- \)harmonic \( (0,q) \) forms with values
in \( L^{k}, \) satisfy

\[
\begin{array}{lr}
\limsup _{k}B_{X}^{q,k}(x)\leq k^{n}B^{q}_{x,\C ^{n}}(0), & \limsup _{k}S_{X}^{q,k}(x)\leq k^{n}S^{q}_{x,\C ^{n}}(0)
\end{array},\]
where 
\[
B^{q}_{x,\C ^{n}}(0)=S^{q}_{x,\C ^{n}}(0)=\frac{1}{\pi ^{n}}1_{X(q)}\left| det_{\omega }(\frac{i}{2}\partial \overline{\partial }\phi )_{x}\right| \]
 and \( \lim _{k}\frac{1}{k^{n}}B_{X}^{q,k}(x)=\lim _{k}\frac{1}{k^{n}}S_{X}^{q,k}(x) \)
if one of the limits exists.
\end{thm}
\begin{proof}
First we will prove that 
\[
\limsup _{k}k^{-n}S_{X}^{q,k}(x)\leq S^{q}_{x,\C ^{n}}(0)\]
 By definition, there is a unit norm sequence \( \alpha _{k} \) in \( \mathcal{H}^{0,q}(X,L^{k}) \)
such that 
\[
\limsup _{k}k^{-n}S_{X}^{q,k}(x)=\limsup _{k}k^{-n}\left| \alpha _{k}(x)\right| ^{2}.\]
Now consider the sequence \( \beta ^{(k)}:=k^{-\frac{1}{2}n}\alpha ^{(k)}, \)
where \( \beta ^{(k)} \) is a form on the ball \( B_{\sqrt{k}R_{k}} \) that
we identify with a form in \( L_{\phi _{0}}^{2}(\C ^{n}), \) by extending with
zero. Note that 
\[
\limsup _{k}\left\Vert \beta ^{(k)}\right\Vert _{\phi _{0}}^{2}=\limsup _{k}k^{-n}\left\Vert \alpha ^{(k)}\right\Vert _{\phi _{0},B_{\sqrt{k}R_{k}}}^{2}\leq \limsup _{k}\left\Vert \alpha _{k}\right\Vert _{X}^{2}=1,\]
 where we have used the norm localization \ref{norm localization}. According
to the previous lemma there is a subsequence of \( \left\{ \beta ^{(k_{j})}\right\}  \)
that converges uniformly with all derivatives to \( \beta  \) on any ball in
\( \C ^{n}, \) where \( \beta  \) is smooth and \( \left\Vert \beta \right\Vert _{\phi _{0}}^{2}\leq 1. \)
Hence, we have that \( \Delta _{\overline{\partial },\phi _{0}}\beta =0, \)
which follows from the expansion \ref{operator localization}, showing that

\[
\limsup _{k}k^{-n}S^{q,k}(x)=\lim _{j}\left| \beta ^{(k_{j})}(0)\right| ^{2}=\left| \beta (0)\right| ^{2}\leq \frac{\left| \beta (0)\right| ^{2}}{\left\Vert \beta \right\Vert ^{2}_{\phi _{0}}}\leq S^{q}_{x,\C ^{n}}(0),\]
Moreover, by proposition \ref{prop:modelB}, 
\[
S^{q}_{x,\C ^{n}}(0)=B^{q}_{x,\C ^{n}}(0)=\frac{1}{\pi ^{n}}1_{X(q)}(x)\left| det_{\omega }(\frac{i}{2}\partial \overline{\partial }\phi )_{x}\right| .\]
 Lemma \ref{lemma: comparing B and S in general}, then shows that \( \begin{array}{lr}
\lim _{k}B_{X}^{q,k}(x)=0
\end{array} \) outside the set \( X(q). \)

Next, if \( x \) is in \( X(q) \) we may assume that \( \lambda _{1} \) up
to \( \lambda _{q} \) are the negative eigenvalues. By \ref{S component in model}
we then have that \( \beta ^{I}=0 \) if \( I\neq (1,2,...,q). \) We deduce
that for \( I\neq (1,2,...,q): \) 
\[
\lim k_{j}^{-n}S_{I}^{k_{j},q}(0)=\lim k_{j}^{-n}\left| \alpha ^{I}_{k_{j}}(0)\right| ^{2}=\left| \beta ^{I}(0)\right| ^{2}=0.\]
This proves that 
\[
\lim k^{-n}S_{I}^{k,q}(0)=0,\]
 if \( I\neq (1,2,...,q). \) Finally, lemma \ref{lemma: comparing B and S in general}
shows that 
\[
\limsup _{k}k^{-n}B_{X}^{k,q}(x)\leq 0+0+...+S^{q}_{x,\C ^{n}}(0)=B^{q}_{x,\C ^{n}}(0),\]
 which finishes the proof of the theorem.
\end{proof}
As a corollary we obtain Demailly's weak holomorphic Morse inequalities:

\begin{cor}
\label{cor:WeakMorse}Suppose that \( X \) is compact. Then 
\[
\dim _{\C }H_{\overline{\partial }}^{0,q}(X,L^{k})\leq k^{n}\frac{(-1)^{q}}{\pi ^{n}}\frac{1}{n!}\int _{X(q)}(\frac{i}{2}\partial \overline{\partial }\phi )^{n}+o(k^{n}).\]

\end{cor}
\begin{proof}
Let us first show that the sequence \( k^{-n}S_{X}^{q,k}(x) \) is dominated
by a constant if \( X \) is compact. Since \( X \) is compact it is enough
to prove this for a sufficintly small neighbourhood of a fixed point \( x_{0}. \)
Now, for a given form \( \alpha _{k} \) in \( \mathcal{H}^{0,q}(X,L^{k}) \)
consider its restriction to a ball of radius \( \frac{1}{\sqrt{k}} \) with
center in \( x_{0}. \) Using Gårding's inequality \ref{proofLemmaGardingA}
as in lemma \ref{lemma:bootstrap}, we see that there is a constant \( C(x_{0}), \)
depending continuosly on \( x_{0} \) such that
\[
\left| \alpha _{k}(x_{0})\right| ^{2}\leq C(x_{0})\left\Vert \alpha ^{(k)}\right\Vert _{\phi _{0},B_{1}}^{2}\]
for \( k \) larger than \( k_{0}(x_{0}), \) say. Morever, we may assume that
the same \( k_{0}(x_{0}) \) works for all \( x \) sufficiently close to \( x_{0}. \)
Using the norm localization \ref{norm localization} we get that 
\[
k_{^{^{-n}}}\left| \alpha _{k}(x_{0})\right| ^{2}\leq 2C(x_{0})\left\Vert \alpha ^{k}\right\Vert _{X}^{2}\]
 for \( k \) larger than \( k_{1}(x_{0}) \) and the same \( k_{1} \) works
for all \( x \) sufficiently close to \( x_{0}. \) This proves that \( k^{-n}S_{X}^{q,k}(x) \)
is dominated by a constant if \( X \) is compact. By \ref{lemma: comparing B and S in general}
and the fact that \( X \) has finite volume, this means that the sequence \( k^{-n}B_{X}^{q,k}(x) \)
is dominated by an \( L^{1}- \)function. Finally, the Hodge theorem shows that
\[
\limsup _{k}\dim _{\C }k^{-n}H_{\overline{\partial }}^{0,q}(X,L^{k})=\limsup _{k}\int _{X}k^{-n}B_{X}^{q,k}\]
 and Fatou's lemma yields, since the sequence \( k^{-n}B_{X}^{q,k} \) is \( L^{1}- \)
dominated, 
\[
\int \limsup k^{-n}B_{X}^{q,k}\leq \frac{(-1)^{q}}{\pi ^{n}}\int _{X(q)}(\frac{i}{2}\partial \overline{\partial }\phi )^{n},\]
where we have used the previous theorem and the fact that 
\[
\frac{1}{n!}(\frac{i}{2}\partial \overline{\partial }\phi )^{n}=(-1)^{q}\left| det_{\omega }(\frac{i}{2}\partial \overline{\partial }\phi )_{x}\right| Vol_{X}\]
 on \( X(q). \)
\end{proof}
\begin{rem}
If there is no point where the curvature two-form is non-degenerate and has
exactly one negative eigenvalue, then Demailly observed (see \cite{d2}) that
combining his weak inequalities with their strong counterparts, gives that 
\[
\dim _{\C }H_{\overline{\partial }}^{0,q}(X,L^{k})=k^{n}\frac{(-1)^{q}}{\pi ^{n}}\frac{1}{n!}\int _{X(0)}(\frac{i}{2}\partial \overline{\partial }\phi )^{n}+o(k^{n}).\]
Now, if one also uses theorem \ref{thm:localweak}\footnote{%
Only the case when \( q=0 \) is needed. This case is considerably more elementary
than the other cases. Indeed, all sections are holomorphic (independently of
\( k) \) and one can use the submean inequality for holomorphic functions,
without invoking the scaled Laplacian.
} one gets the following asymptotic expression for the Bergman kernel function:
\[
\limsup _{k}\frac{1}{k^{n}}B^{0,k}_{X}(x)=\frac{1}{\pi ^{n}}\left| det_{\omega }(\frac{i}{2}\partial \overline{\partial }\phi )_{x}\right| ,\]
almost everywhere on the part of \( X \) where the curvature two-form is strictly
positive (on the complement the limit is zero). This seems to be a new result.
In case \( L \) is stricly positive on all of \( X \) a complete asymptotic
expansion of \( B^{0,k}_{X} \) is known (\cite{z}).
\end{rem}

\section{The model case\label{section:model}}

In this section we will be concerned with \( \C ^{n} \) with its standard metric.
Any smooth function \( \phi  \) defines a hermitian metric on the trivial line
bundle and associated bundles, via \( \left| 1\right| _{\phi }^{2}(z)=e^{-\phi (z)}. \)
Explicitly, this means that if \( \alpha ^{0,q}=\sum _{I}f_{I}d\overline{z}^{I} \)
is a \( (0,q)- \)form on \( \C ^{n} \) , then 
\[
\left| \alpha ^{0,q}\right| _{\phi }^{2}(z)=\sum _{I}\left| f_{I}(z)\right| ^{2}e^{-\phi (z)}.\]
The standard differential operators on smooth functions are extended to operators
on forms, by letting them act componentwise. We denote by \( \frac{\partial }{\partial \overline{z_{i}}}^{*} \)
the formal adjoint of \( \frac{\partial }{\partial \overline{z_{i}}} \) with
respect to the norm induced by \( \phi . \) A partial integration shows that

\begin{equation}
\label{identity: adjoint}
\frac{\partial }{\partial \overline{z_{i}}}^{*}=e^{\phi }(-\frac{\partial }{\partial z_{i}})e^{-\phi }=-\frac{\partial }{\partial z_{i}}+\phi _{z}
\end{equation}
 The following classical commutation relations (\cite{h}) are essential for
what follows: 
\begin{equation}
\label{identity: commutation 1}
\frac{\partial }{\partial \overline{z}_{i}}\frac{\partial }{\partial \overline{z_{j}}}^{*}-\frac{\partial }{\partial \overline{z}_{j}}^{*}\frac{\partial }{\partial \overline{z}_{i}}=\phi _{\overline{z_{i}}z_{j}}
\end{equation}
In this section \( \phi =\phi _{0}:=\sum ^{n}_{i=1}\lambda _{i}\left| z_{i}\right| ^{2}, \)
so that the right hand side simplifies to \( \delta _{ij}\lambda _{ii}. \) 

The next lemma gives an explicit expression for \( \Delta _{\overline{\partial },\phi _{0}}, \)
that will enable us to compute the model Bergman kernel function.

\begin{lem}
We have that 
\[
\Delta _{\overline{\partial }}(fd\overline{z}^{I})=(\sum _{i\in I}\frac{\partial }{\partial \overline{z}_{i}}\frac{\partial }{\partial \overline{z_{i}}}^{*}+\sum _{i\in I^{c}}\frac{\partial }{\partial \overline{z_{i}}}^{*}\frac{\partial }{\partial \overline{z}_{i}})fd\overline{z}^{I}=:(\Delta _{\overline{\partial },I}f)d\overline{z}^{I}),\]
 with respect to the fiber metric \( \phi _{0} \) on the trivial line bundle.
\end{lem}
\begin{proof}
All adjoints and formal adjoints are taken with respect to \( \phi _{0}. \)
The following notation will be used. We let \( d\overline{z}^{i} \) act on
forms by wedge multiplication, and denote the adjoint by \( d\overline{z}^{i*}. \)
Then we have the anti-commutation relations
\begin{equation}
\label{proofofModelLaplaceA}
\overline{dz}^{i}\overline{dz}^{j*}+\overline{dz}^{j*}\overline{dz}^{i}=0\textrm{ if }i\neq j
\end{equation}
Also, 

\begin{equation}
\label{proofModelLaplaceB}
\begin{array}{cc}
(\overline{dz}^{i*}d\overline{z}^{i})d\overline{z}^{J}=\left\{ \begin{array}{rl}
1, & i\notin J\\
0, & i\in J
\end{array}\right.  & \overline{(dz}^{i}\overline{dz}^{i*})d\overline{z}^{J}=\left\{ \begin{array}{rl}
0, & i\notin J\\
1, & i\in J
\end{array}\right. 
\end{array}
\end{equation}
 \( \overline{\partial } \) can now be expressed as
\[
\overline{\partial }=\sum _{i}\frac{\partial }{\partial \overline{z_{i}}}\overline{dz}^{i}=\sum _{i}\overline{dz}^{i}\frac{\partial }{\partial \overline{z_{i}}.}\]
Applying \( * \) to this relation immediately yields
\[
\overline{\partial }^{*}=\sum _{i}\frac{\partial ^{*}}{\partial \overline{z_{i}}}\overline{dz}^{i*}=\sum _{i}\overline{dz}^{i*}\frac{\partial ^{*}}{\partial \overline{z_{i}}.}.\]
Now we compute \( \Delta _{\overline{\partial }}:=\overline{\partial }\overline{\partial }^{*}+\overline{\partial }^{*}\overline{\partial :} \)
\[
\overline{\partial }\overline{\partial }^{*}+\overline{\partial }^{*}\overline{\partial }=\sum _{i,j}\frac{\partial }{\partial \overline{z_{i}}}\frac{\partial ^{*}}{\partial \overline{z_{j}}}\overline{dz}^{i}\overline{dz}^{j*}+\sum _{i,j}\frac{\partial ^{*}}{\partial \overline{z_{j}}}\frac{\partial }{\partial \overline{z_{i}}}\overline{dz}^{j*}\overline{dz}^{i}.\]
 By applying \ref{proofofModelLaplaceA} to second term and splitting the sum
this equals 
\[
\sum _{i}\frac{\partial }{\partial \overline{z_{i}}}\frac{\partial ^{*}}{\partial \overline{z_{i}}}\overline{dz}^{i}\overline{dz}^{i*}+\sum _{i}\frac{\partial ^{*}}{\partial \overline{z_{i}}}\frac{\partial }{\partial \overline{z_{i}}}\overline{dz}^{i*}\overline{dz}^{i}+\sum _{i\neq j}(\frac{\partial }{\partial \overline{z_{i}}}\frac{\partial ^{*}}{\partial \overline{z_{j}}}-\frac{\partial ^{*}}{\partial \overline{z_{j}}}\frac{\partial }{\partial \overline{z_{i}}}\overline{)dz}^{i}\overline{dz}^{j*}.\]
According to the commutation relation \ref{identity: commutation 1} the second
sum vanishes and \ref{proofModelLaplaceB} finally gives that 
\[
(\overline{\partial }\overline{\partial }^{*}+\overline{\partial }^{*}\overline{\partial })(fd\overline{z}^{I})=(\sum _{i\in I}\frac{\partial }{\partial \overline{z}_{i}}\frac{\partial }{\partial \overline{z_{i}}}^{*\phi _{0}}+\sum _{i\in I^{c}}\frac{\partial }{\partial \overline{z_{i}}}^{*\phi _{0}}\frac{\partial }{\partial \overline{z}_{i}})fd\overline{z}^{I}.\]
 
\end{proof}
In order to be able to integrate partially in \( \C ^{n} \) without getting
boundary terms, the following lemma is useful. See \cite{gro} for a slightly
more general proof. The point is to choose \( \chi _{R}(z)=(\chi (\frac{z}{R}))^{2}, \)
where \( \chi  \) is a smooth compact supported function, that equals \( 1 \)
on the unit ball, say.

\begin{lem}
\label{lemma: gromov-trick}There is an exhaustion sequence \( \chi _{R} \)
such that for any smooth function \( f \) and smooth \( (0,q)- \) form \( \alpha  \)
with \( f, \) \( \frac{\partial }{\partial \overline{z_{i}}}f, \) \( \alpha  \)
and \( (\overline{\partial }+\overline{\partial }^{*})\alpha  \) in \( L^{2}_{\phi _{0}}(\C ^{n}) \)
\[
\lim _{R}\left\langle \Delta _{\overline{\partial }}\alpha ,\chi _{R}\alpha \right\rangle =\left\Vert (\overline{\partial }+\overline{\partial }^{*})\alpha \right\Vert _{\phi _{0}}^{2},\, \lim _{R}\left\langle \frac{\partial }{\partial \overline{z}_{i}}^{*}\frac{\partial }{\partial \overline{z_{i}}}f,\chi _{R}f\right\rangle =\left\Vert \frac{\partial }{\partial \overline{z_{i}}}f\right\Vert ^{2}\]
and similarly for \( \left\Vert \frac{\partial ^{*}}{\partial \overline{z_{i}}}f\right\Vert ^{2}. \)
\end{lem}
Next, we turn to the Bergman kernel function \( B^{q}_{\phi _{0},\C ^{n}}(z) \)
and the extremal function \( S_{\phi _{0},\C ^{n}}^{q}(z) \). They are defined
as in section \ref{section: dimension, ...} but with respect to the space \( \mathcal{H}_{\phi _{0}}^{0,q}(\C ^{n}) \)
of all \( (0,q)- \)forms in \( \C ^{n}, \) that are \( \overline{\partial }- \)harmonic
with respect to the fiber metric \( \phi _{0} \) and have finite \( L^{2} \)
norm with respect to \( \phi _{0}. \) Note that with notation as in section
\ref{subsection statement of main result}, 
\[
B^{q}_{\phi _{0},\C ^{n}}(z)=B^{q}_{x,\C ^{n}}(z)\]
 if \( \phi (z)=\phi _{0}(z)+O(\left| z\right| ^{3}) \) is the local expression
of the fiber metric of \( L \) at the point \( x \) in \( X. \) 

The proof of the following proposition is based on a reduction to the case when
\( q=0 \) and we are considering holomorphic functions in the so called Fock
space. Then it is well-known that 
\begin{equation}
\label{bergmankernel for holomorphic f}
B^{q}_{x,\C ^{n}}(z)=S_{\C ^{n}}^{q}(0)=\frac{1}{\pi ^{n}}\left| \lambda _{1}\right| \left| \lambda _{2}\right| \cdot \cdot \cdot \left| \lambda _{n}\right| .
\end{equation}
 Indeed, if \( f \) is holomorphic, then \( \left| f\right| ^{2} \) is subharmonic.
Hence, 
\begin{equation}
\label{ineq.: submean ineq.1}
\int _{\Delta _{R}}\left| f(0)\right| ^{2}e^{-k\phi _{0}(z)}\leq \int _{\Delta _{R}}\left| f(z)\right| ^{2}e^{-k\phi _{0}(z)},
\end{equation}
 where \( \Delta _{R} \) is a polydisc of radius \( R \) and where we have
used that \( \phi _{0} \) is radial in each variable. Letting \( R \) tend
to infinity, shows \ref{bergmankernel for holomorphic f} for \( S_{\C ^{n}}^{q}(0). \)

\begin{prop}
\label{prop:modelB}Assume that \( q \) of the numbers \( \lambda _{i} \)
are negative and the rest are positive. Then
\[
B^{q}_{\phi _{0},\C ^{n}}(0)=S_{\phi _{0},\C ^{n}}^{q}(0)=\frac{1}{\pi ^{n}}\left| \lambda _{1}\right| \left| \lambda _{2}\right| \cdot \cdot \cdot \left| \lambda _{n}\right| ,\]
Otherwise there are no \( \overline{\partial }- \)harmonic \( (0,q)- \)forms
in \( L_{\phi _{0}}^{2}(\C ^{n}), \) i.e.
\[
B^{q}_{\C ^{n}}(0)=S_{\C ^{n}}^{q}(0)=0.\]
With notation as in section \ref{subsection statement of main result} this
means that 
\[
B^{q}_{x,\C ^{n}}(0)=S^{q}_{x,\C ^{n}}(0)=\frac{1}{\pi ^{n}}1_{X(q)}(x)\left| det_{\omega }(\frac{i}{2}\partial \overline{\partial }\phi )_{x}\right| \]

\end{prop}
\begin{proof}
Suppose that \( \alpha ^{0,q}=\sum _{I}f_{I}d\overline{z}^{I},\alpha ^{0,q}\in L_{\phi _{0}}^{2}(\C ^{n}) \)
and \( \Delta _{\overline{\partial ,}\phi _{0}}\alpha ^{0,q}=0. \) Then \( \left\langle \Delta _{\overline{\partial }}\alpha ^{0,q},\chi _{R}\alpha ^{0,q}\right\rangle _{\phi _{0}}=0 \)
and we get that 
\[
\sum _{I}\left\langle \Delta _{I}f_{I},\chi _{R}f_{I}\right\rangle _{\phi _{0}}=0.\]
Letting \( R\rightarrow \infty  \) and using lemma \ref{lemma: gromov-trick}
shows that
\begin{equation}
\label{proodmodela}
\begin{array}{cc}
\left\Vert \frac{\partial ^{*\phi _{0}}}{\partial \overline{z_{i}}}f_{I}\right\Vert _{\phi _{0}}^{2}=0,\, i\in I & \left\Vert \frac{\partial }{\partial \overline{z_{i}}}f_{I}\right\Vert _{\phi _{0}}^{2}=0,\, i\in I^{c}
\end{array}.
\end{equation}
 Let \( F_{I}(\zeta ):=e^{-\sum _{i\in I}\lambda _{i}\left| z_{i}\right| ^{2}}f_{I}, \)
where \( \zeta _{i}=\overline{z}_{i} \) if \( i\in I \) and \( \zeta _{i}=z_{i} \)
if \( i\in I^{c}. \) Then \ref{proodmodela} and \ref{identity: adjoint} implies
that \( F_{I} \) is holomorphic: 
\[
\frac{\partial }{\partial \overline{z}_{i}}F_{I}=0\textrm{ for all }i.\]
 Moreover,
\begin{equation}
\label{proofmodelb}
\left| f_{I}\right| ^{2}_{\phi _{0}}=\left| F_{I}\right| _{\Phi _{I}}^{2},\, \Phi _{I}(z):=\sum _{i\in I}-\lambda _{i}\left| z_{i}\right| ^{2}+\sum _{i\in I^{c}}\lambda _{i}\left| z_{i}\right| ^{2}.
\end{equation}
 Now assume that it is not the case that \( q \) of the numbers \( \lambda _{i} \)
are negative and the rest are positive. Then \( \Phi _{I}(z)=\sum ^{n}_{i=1}\Lambda ^{I}_{i}\left| z_{i}\right| ^{2}, \)
with some \( \Lambda ^{I}_{i}:=\Lambda ^{I}_{i_{I}}\leq 0. \) By assumption,
\( \int \left| f_{I}\right| ^{2}e^{-\phi }<\infty . \) Now \ref{proofmodelb}
and Fubini-Tonelli's theorem give that 
\[
\int \left| F_{I}(0,..,z_{i_{I}},0,..)\right| _{\phi _{I}}^{2}e^{-\Lambda ^{I}_{i_{I}}\left| z_{i_{I}}\right| ^{2}}<\infty .\]
Since \( \frac{\partial }{\partial \overline{z_{i_{I}}}}F_{I}=0 \) and \( \Lambda ^{I}_{i_{I}}\leq 0 \)
this forces \( F_{I}(0,..,z_{i_{I}},0,..)\equiv 0 \) and in particular \( f_{I}(0)=0, \)\footnote{%
this follows for example from the submean inequality \ref{ineq.: submean ineq.1}
}which proves that \( B^{q}_{\C ^{n}}(0)=S_{\C ^{n}}^{q}(0)=0. \) 

Finally, assume that \( q \) of the numbers \( \lambda _{i} \) are negative
and the rest are positive. We may assume that \( \lambda _{1} \) up to \( \lambda _{q} \)
are the negative ones. Then the same argument as the one above gives that \( f_{I}(0)=0 \)
if \( I\neq I_{0}:=(1,2,...,q). \) Now since \( \Lambda ^{I_{0}}_{i}>0 \)
for all \( i \) we get

\[
S_{\C ^{n}}^{q}(0)=\textrm{sup}\frac{\left| \alpha (x)\right| ^{2}}{\left\Vert \alpha \right\Vert _{\phi _{0}}^{2}}=\textrm{sup}\frac{\left| f_{I_{0}}(x)\right| ^{2}}{\left\Vert f_{I_{0}}\right\Vert _{\phi _{0}}^{2}}=\textrm{sup}\frac{\left| F_{I_{0}}(x)\right| ^{2}}{\left\Vert F_{I_{0}}\right\Vert _{\Phi _{I}}^{2}}=\frac{\Lambda _{1}\Lambda _{2}\cdot \cdot \cdot \Lambda _{n}}{\pi ^{n}},\]
 where we have used \ref{bergmankernel for holomorphic f} in the last step.
Since \( \Lambda _{i}=\left| \lambda _{i}\right|  \)this proves the statement
about \( S_{\C ^{n}}^{q}(0). \) The proof is finished by observing that \( S_{\C ^{n}}^{q}(0)=B^{q}_{\C ^{n}}(0). \)
This follows from lemma \ref{lemma: comparing B and S in general}, but it is
also easy to see directly in this special case, since all the components \( F_{I} \)
vanish if \( I\neq (1,2,...,q). \)
\end{proof}
\begin{rem}
The statement \ref{S component in model} also follows from the previous proof.

A similar argument to the one in the previous proof shows that if \( \partial \overline{\partial }\phi _{0} \)
is non-degenerate and if it is \emph{not} the case that \( q \) of the numbers
\( \lambda _{i} \) are negative and \( n-q \) of the numbers are positive,
then there is an apriori estimate of the form 
\[
\left\Vert \alpha ^{q}\right\Vert ^{2}\leq C_{q}(\left\Vert \overline{\partial }\alpha ^{q}\right\Vert ^{2}+\left\Vert \overline{\partial }^{*}\alpha ^{q}\right\Vert ^{2}),\]
 where the norms are taken with respect to \( \phi _{0}. \) The result appears
already in Hörmander's seminal paper \cite{h} and it can be used to give a
direct proof of the fact that the global Bergman kernel function \( B^{q,k}_{X} \)
vanishes (modulo terms of order \( o(k^{n}) \) at a point outside \( X(q), \)
where the curvature two-form is non-degenerate. 
\end{rem}

\section{The strong holomorphic morse inequalities}

In this section \( X \) is assumed to be compact. While the weak holomorphic
Morse inequalities give estimates on individual cohomology groups, their strong
counter parts say that 
\begin{equation}
\label{strong inequalities}
\sum ^{q}_{j=0}(-1)^{q-j}\dim _{\C }H^{j}(X,L^{k})\leq k^{n}\frac{(-1)^{q}}{\pi ^{n}}\frac{1}{n!}\int _{X(\leq q)}(\frac{i}{2}\partial \overline{\partial }\phi )^{n}+o(k^{n}).
\end{equation}
Let \( \mathcal{H}_{\leq \nu _{k}}^{q}(X,L^{k}) \) denote the space spanned
by the eigenforms of \( \Delta _{\overline{\partial }} \) whose eigenvalues
are bounded by \( \nu _{k} \) and denote by \( B_{\leq \nu _{k}}^{q,k} \)
the Bergman kernel function of the space \( \mathcal{H}_{\leq \nu _{k}}^{q}(X,L^{k}). \)
Extending the previous methods (section \ref{section:weakMorse}) we will show
the asymptotic equality
\[
\lim _{k}k^{-n}B_{\leq \nu _{k}}^{q,k}(x)=\frac{1}{\pi ^{n}}1_{X(q)}\left| det_{\omega }(\frac{i}{2}\partial \overline{\partial }\phi )_{x}\right| ,\]
 where \( \left\{ \nu _{k}\right\}  \) is a judicially chosen sequence. From
its integrated version 
\[
\dim _{\C }\mathcal{H}_{\leq \nu _{k}}^{q}(X,L^{k})=k^{n}\frac{(-1)^{q}}{\pi ^{n}}\frac{1}{n!}\int _{X(q)}(\frac{i}{2}\partial \overline{\partial }\phi )^{n}+o(k^{n})\]
 one can then deduce the strong inequalities \ref{strong inequalities} as in
\cite{d1}. The basic idea is as follows. First we obtain an \emph{upper} bound
on the corresponding Bergman kernel function, by a direct generalizaton of the
harmonic case: the local weak Morse inequalities. We just have to make sure
that terms of the form \( \left\Vert (\Delta _{\overline{\partial }}^{(k)})^{m}\beta ^{(k)}\right\Vert _{\sqrt{k}R_{k}}^{2}, \)
which are now non-zero, tend to zero with \( k. \) Recall, that \( \Delta _{\overline{\partial }}^{(k)} \)
denotes the \( \overline{\partial } \)-Laplacian with respect to the scaled
metrics on the ball \( B_{\sqrt{k}R_{k}}. \) In fact, for scaling reasons they
give contributions which are polynomial in\( \frac{\nu _{k}}{k}. \) This dictates
the choice \( \nu _{k}=\mu _{k}k \) with \( \mu _{k} \) tending to zero with
\( k. \) The last step is to get a \emph{lower} bound on the spectral density
function on \( X(q), \) which amounts to proving the existence of a unit norm
sequence \( \left\{ \alpha _{k}\right\}  \) in \( \mathcal{H}_{\leq \nu _{k}}^{i}(X,L^{k}) \)
with 
\[
\left| a_{k}(x)\right| ^{2}=k^{n}\frac{1}{\pi ^{n}}\left| det_{\omega }(\frac{i}{2}\partial \overline{\partial }\phi )_{x}\right| +o(k^{n}).\]
To this end we first take a sequence \( \left\{ \alpha _{k}\right\}  \) of
\( (0,q)- \)forms on \( \C ^{n}, \) which are harmonic with respect to the
\emph{flat} metric and the fiber metric \( k\phi _{0}, \) having the corresponding
property. The point is that the mass of these forms concentrate around \( 0, \)
when \( k \) tends to infinity. Hence, by cutting down their support to small
decreasing balls we obtain global forms \( \widetilde{\alpha _{k}} \) on \( X, \)
that are of unit norm in the limit. Moreover, their Laplacians are ``small''.
By cutting the high frequencies, that is projecting on \( \mathcal{H}_{\leq \nu _{k}}^{i}(X,L^{k}), \)
we finally get the sought after sequence.

We now proceed to carry out the details of the argument sketched above.

\begin{prop}
\label{prop:Bspectral<}Assume that \( \mu _{k}\rightarrow 0. \) Then the following
estimate holds:
\[
B_{\leq \mu _{k}k}^{q,k}(x)\leq k^{n}\frac{1}{\pi ^{n}}1_{X(q)}\left| det_{\omega }(\frac{i}{2}\partial \overline{\partial }\phi )_{x}\right| +o(k^{n}).\]
 
\end{prop}
\begin{proof}
The proof is a simple modification of the proof of the local holomorphic Morse
inequalities and in what follows these modifications will be presented. The
difference is that \( \alpha _{k} \) is in \( \mathcal{H}_{\leq \mu _{k}k}^{0,q}(X,L^{k}) \)
and we have to make sure that the all terms of the form \( (\Delta _{\overline{\partial }}^{(k)})^{m}(k^{-\frac{n}{2}}\alpha ^{(k)})=(\Delta _{\overline{\partial }}^{(k)})^{m}(\beta ^{(k)}) \)
vanish in the limit. But for any ball \( B, \)
\[
\left\Vert (\Delta _{\overline{\partial }}^{(k)})^{m}\beta ^{(k_{j})}\right\Vert _{\phi _{0},B}^{2}\leq k^{-n}\left\Vert (\Delta _{\overline{\partial }}^{(k)})^{m}\alpha ^{(k_{j})}\right\Vert _{\phi _{0},B_{\sqrt{k}R_{k}}}^{2}\lesssim \leq k^{-2m}\left\Vert (\Delta _{\overline{\partial }})^{m}\alpha _{k}\right\Vert _{X}^{2}\]
and the last term is bounded by a sequence tending to zero:
\[
k^{-2m}(\mu _{k}k)^{2m}\rightarrow 0,\]
since by assumption \( \alpha _{k} \) is of unit norm and in \( \mathcal{H}_{\leq \mu _{k}k}^{0,q}(X,L^{k}) \)
and \( \mu _{k}\rightarrow 0. \) Gårding's inequality as in \ref{proofLemmaGardingA}
gives that

\[
\left\Vert \beta ^{(k_{j})}\right\Vert ^{2}_{\phi _{0},B,2,m}\leq C\left( \left\Vert \beta ^{(k_{j})}\right\Vert ^{2}_{\phi _{0},B}+\left\Vert (\Delta _{\overline{\partial }}^{(k)})^{m}\beta ^{(k_{j})}\right\Vert ^{2}_{\phi _{0},B}\right) \lesssim (C+(\mu _{k_{j}})^{2m})\lesssim C,\]
which shows that the conclusion of lemma \ref{lemma:bootstrap} is still valid.
Finally, \( \Delta _{\overline{\partial }}\beta =0 \) as before and the rest
of the argument goes through word by word.
\end{proof}
The next lemma provides the sequence that takes the right values at a given
point \( x \) in \( X(q), \) with ``small'' Laplacian, that was referred
to in the beginning of the section.

\begin{lem}
\label{lemma:lowenergysequence}Let \( c_{\phi }(x):=\frac{1}{\pi ^{n}}\left| det_{\omega }(\frac{i}{2}\partial \overline{\partial }\phi )_{x}\right| . \)
For any point \( x_{0} \) in \( X(q) \) there is a sequence \( \left\{ \alpha _{k}\right\}  \)
such that \( a_{k} \) is in \( \Omega ^{0,q}(X,L^{k}) \) with 
\[
\begin{array}{ll}
(i) & \left| a_{k}(x)\right| ^{2}=k^{n}c_{\phi }(x)\\
(ii) & \lim _{k}\left\Vert \alpha _{k}\right\Vert ^{2}=1\\
(iii) & \left\Vert k^{-m}(\Delta _{\overline{\partial }})^{m}\alpha _{k}\right\Vert ^{2}_{X}=0
\end{array}\]
 Moreover, there is a sequence \( \delta _{k} \) , independent of \( x_{0} \)
and tending to zero, such that 
\[
(iv)\, \left\langle k^{-1}\Delta _{\overline{\partial }}\alpha _{k},\alpha _{k}\right\rangle _{X}\leq \delta _{k}\]

\end{lem}
\begin{proof}
We may assume that the first \( q \) eigenvalues \( \lambda _{i,x_{0}} \)
are negative, while the remaining eigenvalues are positive\( . \) Define the
following form in \( \C ^{n}: \)
\[
\beta (w)=\left( \frac{^{\left| \lambda _{1}\right| \left| \lambda _{2}\right| \cdot \cdot \cdot \left| \lambda _{n}\right| }}{\pi ^{n}}\right) ^{\frac{1}{2}}e^{+\sum ^{q}_{i=1}\lambda _{i}\left| w_{i}\right| ^{2}}d\overline{w_{1}}\wedge d\overline{w_{2}}\wedge ...\wedge d\overline{w_{q}},\]
 so that\footnote{%
compare the proof of proposition \ref{prop:modelB}.
} \( \left| \beta \right| ^{2}_{\phi _{0}}=\frac{^{\left| \lambda _{1}\right| \left| \lambda _{2}\right| \cdot \cdot \cdot \left| \lambda _{n}\right| }}{\pi ^{n}}e^{-\sum ^{n}_{i=1}\left| \lambda _{i}\right| \left| w_{i}\right| ^{2}} \)
and \( \left\Vert \beta \right\Vert _{\phi _{0},\C ^{n}}=1. \) Observe that
\( \beta  \) is in \( L_{\phi _{0}}^{2,m}, \) the Sobolev space with \( m \)
derivatives in \( L_{\phi _{0}}^{2}, \) for all \( m \) Now define \( \alpha _{k} \)
on \( X \) by
\[
\alpha _{k}(z):=k^{\frac{n}{2}}\chi _{k}(\sqrt{k}z)\beta (\sqrt{k}z),\]
 where \( \chi _{_{k}}(w)=\chi (\frac{w}{\sqrt{k}R_{k}}) \) and \( \chi  \)
is a a smooth function supported on the unit ball, which equals one on the ball
of radius \( \frac{1}{2}. \) Thus \( \left| a_{k}(x_{0})\right| ^{2}=k^{n}c_{\phi }(x), \)
showing \( (i). \) To see \( (ii) \) note that 
\begin{equation}
\label{tail}
\left\Vert \alpha _{k}\right\Vert _{X}^{2}=_{k}\left\Vert \chi _{k}\beta \right\Vert ^{2}_{\phi _{0},\C ^{n}}=\left\Vert \beta \right\Vert ^{2}_{\phi _{0},\frac{1}{2}\sqrt{k}R_{k}}+\left\Vert \chi _{k}\beta \right\Vert ^{2}_{\phi _{0},\geq \frac{1}{2}\sqrt{k}R_{k}}
\end{equation}
 and the ``tail'' \( \left\Vert \beta \right\Vert ^{2}_{\phi _{0},\geq \frac{1}{2}\sqrt{k}R_{k}} \)
tends to zero, since \( \beta  \) is in \( L_{\phi _{0}}^{2}(\C ^{n}) \) and
\( \sqrt{k}R_{k} \) tends to infinity. 

Next, we show \( (iii). \) Changing variables (proposition \ref{scaling transf for laplace})
and using \ref{operator localization} gives 
\begin{equation}
\label{power of laplace of alpha}
\left\Vert k^{-m}(\Delta _{\overline{\partial }})^{m}\alpha _{k}\right\Vert ^{2}_{X}\lesssim \left\Vert (\Delta _{\overline{\partial }}^{(k)})^{m-1}(\Delta _{\overline{\partial },\phi _{0}}+\epsilon _{k}\mathcal{D}_{k})(\chi _{k}\beta ))\right\Vert _{\phi _{0},\sqrt{k}R_{k}}^{2},
\end{equation}
where \( \mathcal{D}_{k} \) is a second order \( PDO \), whose coefficients
have derivatives that are uniformly bounded in \( k. \) To see that this tends
to zero first observe that 
\begin{equation}
\label{Laplace of beta}
\left\Vert (\Delta _{\overline{\partial }}^{(k)})^{m-1}(\Delta _{\overline{\partial },\phi _{0}}\chi _{k}\beta )\right\Vert _{\phi _{0},\sqrt{k}R_{k}}^{2}
\end{equation}
tends to zero. Indeed, \( \beta  \) has been chosen so that \( \Delta _{\overline{\partial },\phi _{0}}\beta =0. \)
Moreover, \( \Delta _{\overline{\partial },\phi _{0}} \) is the square of the
first order operator \( \overline{\partial }+\overline{\partial }^{*,\phi _{0}}, \)
which also annihilates \( \beta  \) and obeys a Leibniz like rule, showing
that 
\[
\Delta _{\overline{\partial },\phi _{0}}\chi _{k}\beta =\gamma _{k}\beta \]
 where \( \gamma _{k} \) is a form, uniformly bounded in \( k, \) and supported
outside the ball \( B_{\frac{1}{2}\sqrt{k}R_{k}} \) (\( \gamma _{k} \) contains
second derivatives of \( \chi _{k}). \) Now using \ref{operator localization}
again we see that \ref{Laplace of beta} is bounded by the norm of 
\[
\gamma _{k}p(w,\overline{w})\beta ,\]
 where \( p \) is a polynomial. Thus \ref{power of laplace of alpha} can be
estimated by the ``tail'' of a convergent integral, as in \ref{tail} - the
polynomial does not effect the convergence - which shows that \ref{Laplace of beta}
tends to zero. To finish the proof of \( (iii) \) it is now enough to show
that
\[
\left\Vert (\Delta _{\overline{\partial }}^{(k)})^{m-1}\mathcal{D}_{k}(\chi _{k}\beta )\right\Vert _{\phi _{0},\sqrt{k}R_{k}}^{2},\]
 is uniformly bounded. As above one sees that the integrand is bounded by the
norm of 
\[
q(w,\overline{w})\beta ,\]
 for some polynomial \( q, \) which is finite as above. 

To prove \( (iv) \) observe that, as above, 
\[
\left\langle k^{-1}\Delta _{\overline{\partial }}\alpha _{k},\alpha _{k}\right\rangle =\left\Vert \frac{1}{\sqrt{k}}(\overline{\partial }+\overline{\partial }^{*})\alpha _{k}\right\Vert ^{2}_{X}\sim \left\Vert (\overline{\partial }+\overline{\partial }^{*(k)})(\chi _{\sqrt{k}R_{k}}\beta )\right\Vert ^{2}_{\sqrt{k}R_{k}}.\]
 Hence, by Leibniz' rule 
\[
\left\langle k^{-1}\Delta _{\overline{\partial }}\alpha _{k},\alpha _{k}\right\rangle \lesssim \left\Vert (\chi _{\sqrt{k}R_{k}}(\overline{\partial }+\overline{\partial }^{*(k)})\beta \right\Vert ^{2}_{\sqrt{k}R_{k}}+\frac{C}{(\sqrt{k}R_{k})^{2}}\left\Vert \beta \right\Vert ^{2}_{\C ^{n}}.\]
Clearly, there is an expansion for the first order operator \( (\overline{\partial }+\overline{\partial }^{*(k)}) \)
as in \ref{operator localization}, giving
\[
\left\langle k^{-1}\Delta _{\overline{\partial }}\alpha _{k},\alpha _{k}\right\rangle \lesssim \varepsilon _{k}\left( \left\Vert \beta \right\Vert ^{2}+\sum ^{2n}_{i=1}\left\Vert \partial _{i}\beta \right\Vert ^{2}\right) +\frac{C}{(\sqrt{k}R_{k})^{2}}\left\Vert \beta \right\Vert ^{2}.\]
 Note that even if \( \left\Vert \beta \right\Vert ^{2} \) is independent of
the eigenvalues \( \lambda _{i,x_{0}}, \) the norms \( \left\Vert \partial _{i}\beta \right\Vert ^{2} \)do
depend on the eigenvalues, and hence on the point \( x_{0}. \) But the dependence
amounts to a factor of eigenvalues and since \( X \) is compact, we deduce
that \( \left\Vert \partial _{i}\beta \right\Vert ^{2} \)is bounded by a constant
independent of the point \( x_{0}. \) This shows that \( \left\langle k^{-1}\Delta _{\overline{\partial }}\alpha _{k},\alpha _{k}\right\rangle _{X}\leq \delta _{k}. \)
Note that \( \epsilon _{k} \) also can be taken to be independent of the point
\( x_{0}, \) by a similar argument.
\end{proof}
\begin{prop}
\label{prop:Bspectral>}Assume that the sequence \( \mu _{k} \) is such that
\( \mu _{k}\neq 0 \) and \( \frac{\delta _{k}}{\mu _{k}}\rightarrow 0, \)
where \( \delta _{k} \) is the sequence appearing in lemma \ref{lemma:lowenergysequence}.
Then, for any point \( x \) in \( X(q), \) the following holds 
\[
\liminf k^{-n}B_{\leq \mu _{k}k}^{q,k}(x)\geq \frac{1}{\pi ^{n}}\left| det_{\omega }(\frac{i}{2}\partial \overline{\partial }\phi )_{x}\right| \]

\end{prop}
\begin{proof}
Let \( \left\{ \alpha _{k}\right\}  \) be the sequence that lemma \ref{lemma:lowenergysequence}
provides and decompose it with respect to the orthogonal decomposition \( \Omega ^{0,q}(X,L^{k})=\mathcal{H}_{\leq \mu _{k}k}^{q}(X,L^{k})\oplus \mathcal{H}_{>\mu _{k}k}^{q}(X,L^{k}), \)
induced by the spectral decomposition of the elliptic operator \( \Delta _{\overline{\partial }}: \)
\[
\alpha _{k}=\alpha _{1,k}+\alpha _{2,k}\]
First, we prove that 
\begin{equation}
\label{alphatwo vanishes}
\lim _{k}k^{-n}\left| \alpha _{2}^{(k)}(0)\right| ^{2}=0.
\end{equation}
As in the proof of lemma \ref{lemma:bootstrap} we have that
\[
k^{-n}\left| \alpha _{2}^{(k)}(0)\right| ^{2}\leq C(x)\left( k^{-n}\left\Vert \alpha _{2}^{(k)}\right\Vert _{B_{1}}^{2}+k^{-n}\left\Vert (\Delta _{\overline{\partial }}^{(k)})^{m}\alpha _{2}^{(k)}\right\Vert _{,B_{1}}^{2}\right) .\]
To see that the first term tends to zero, observe that by the spectral decomposition
of \( \Delta _{\overline{\partial }}: \)
\[
\left\Vert \alpha _{2,k}\right\Vert _{X}^{2}\leq \frac{1}{\mu _{k}k}\left\langle \Delta _{\overline{\partial }}\alpha _{2,k},\alpha _{2,k}\right\rangle _{X}\leq \frac{1}{\mu _{k}}\left\langle k^{-1}\Delta _{\overline{\partial }}\alpha _{k},\alpha _{k}\right\rangle _{X}\leq \frac{\delta _{k}}{\mu _{k}}\]
Furthermore, the second term also tends to zero: 
\[
k^{-n}\left\Vert (\Delta _{\overline{\partial }}^{(k)})^{m}\alpha _{2}^{(k)}\right\Vert _{B_{1}}^{2}\leq \left\Vert k^{-m}(\Delta _{\overline{\partial }})^{m}\alpha _{2,k}\right\Vert _{X}^{2}\leq \left\Vert k^{-m}(\Delta _{\overline{\partial }})^{m}\alpha _{k}\right\Vert _{X}^{2}\rightarrow 0.\]
 by \( (iii) \) in lemma \ref{lemma:lowenergysequence}. Finally, 
\[
k^{-n}S^{q,k}_{\mu _{k}k}(x)\geq k^{-n}\frac{\left| \alpha _{k,1}(0)\right| ^{2}}{\left\Vert \alpha _{k,1}\right\Vert _{X}^{2}}\geq k^{-n}\left| \alpha _{1,k}(0)\right| ^{2}=k^{-n}\left| \alpha _{k}(0)-\alpha _{2,k}(0)\right| ^{2}.\]
By \ref{alphatwo vanishes} this tends to the limit of \( k^{-n}\left| \alpha _{k}(0)\right| ^{2}, \)
which proves the proposition according to \( (i) \) in lemma \ref{lemma:lowenergysequence}
and lemma \ref{lemma: comparing B and S in general}.
\end{proof}
Now we can prove the following asymptotic equality:

\begin{thm}
\label{thm:pointwisestrongMorse} Let \( (X,\omega ) \) be a compact hermitian
manifold. Then 
\[
\lim k^{-n}B_{\leq \mu _{k}k}^{q,k}(x)=\frac{1}{\pi ^{n}}1_{X(q)}\left| det_{\omega }(\frac{i}{2}\partial \overline{\partial }\phi )_{x}\right| ,\]
 for some sequence \( \mu _{k} \) tending to zero. 
\end{thm}
\begin{proof}
Let \( \mu _{k}:=\sqrt{\delta _{k}.} \) The theorem then follows immediately
from proposition \ref{prop:Bspectral<} and proposition \ref{prop:Bspectral>}
if \( x \) is in \( X(q). \) If \( x \) is outside of \( X(q) \) then the
upper bound given by \ref{prop:Bspectral<} shows that \( \lim _{k}B_{\leq \mu _{k}k}^{q,k}(x)=0, \)
which finishes the proof of the theorem.
\end{proof}
\begin{acknowledgement}
The author whishes to thank his advisor Bo Berndtsson for stimulating and enlightening
discussions and for his positive attitude. Furthermore, the author is grateful
to Johannes Sjöstrand for comments on an early draft of the manuscript.
\end{acknowledgement}

\end{document}